\documentclass[12pt]{amsart}

\usepackage{latexsym}
\usepackage{amsmath}
\usepackage{amssymb}
\usepackage{amscd}
\usepackage{psfig}
\usepackage{multicol}
\usepackage{amsfonts}

\newtheorem{theorem}{Theorem}[section]

\newtheorem{lemma}{Lemma}[section]

\def\R{\mathbb{R}}
\def\C{\mathbb{C}}
\def\iso{\cong}

\def\t{\mathfrak{t}}
\def\g{\mathfrak{g}}

\title[Cohomology of Reductions]{An effective algorithm for the
cohomology ring of symplectic reductions}\thanks{This work was partially supported by a
National Science Foundation Postdoctoral Fellowship}

\author{R. F. Goldin\\
{\em G\lowercase{eorge} M\lowercase{ason} U\lowercase{niversity}}}

\begin{document}


\begin{abstract} Let $G$ be a compact torus acting on a compact
symplectic manifold $M$ in a Hamiltonian fashion, and $T$ a subtorus
of $G$. We
prove that the kernel of $\kappa:H_G^*(M)\rightarrow
H^*(M/\!/G)$ is generated by a small number of classes $\alpha\in
H_G^*(M)$ satisfying very explicit restriction properties. Our main tool is
the {\em equivariant Kirwan map}, a natural map from
the $G$-equivariant cohomology of $M$ to the $G/T$-equivariant
cohomology of the symplectic reduction of $M$ by $T$. We show this map is
surjective. This is an equivariant version of the well-known result
that the (nonequivariant) Kirwan map $\kappa:H_G^*(M)\rightarrow
H^*(M/\!/G)$ is surjective. We also compute the kernel of the
equivariant Kirwan map, generalizing the result due to Tolman and
Weitsman \cite{TW:abelianquotient} in the case $T=G$ and allowing us to apply their methods inductively. This result is new even in the case that $\dim T = 1$.
We close
with a worked example: the cohomology ring of the product of two $\C P^2$s,
quotiented by the diagonal 2-torus action.
\end{abstract}

\maketitle

\section{Introduction and Statement of Results}\label{se:intro}

Let $(M,\omega)$ be a symplectic manifold with an
action of a compact torus $G$. A {\em moment map} is an invariant
map
$$
\Phi: M \longrightarrow \g^*
$$
which intertwines the group action and the symplectic form by the {\em
moment map condition}
\begin{equation}\label{eq:momentmapcondition}
\omega(\cdot, X_\xi) = d\langle\Phi,\xi\rangle,
\end{equation}
where $\xi \in \g$, $X_\xi$ is the vector field on $M$ generated by
$\xi$, and $\langle\ ,\ \rangle$ is the pairing of $\g^*$ with
$\g$. We also write $\Phi^\xi$ to indicate $\langle\Phi,\xi\rangle$,
the $\xi$-component of $\Phi$. Condition~(\ref{eq:momentmapcondition})
and the nondegeneracy of $\omega$ imply that singular points of $\Phi$
occur when $X_\xi=0$, or when a subtorus of G acts trivially. If
$C\subset M$ is a component of the fixed point set of $G$ on $M$, then
$\Phi(C)$ has constant value in $\g^*$.

For $M$ a compact manifold, the image of $\Phi$ is a convex polytope
in $\g^*$ \cite{GS:convexity}. Let $d=\dim G$. For any subtorus
$T\subset G$, denote by $M^{T}$ the fixed point set of $M$ under
$T$. We say $T$ is {\em generic} if $M^T = M^G$. The terminology is
appropriate as there are a finite number of subtori $T$ such that
$M^T$ contains but does not equal $M^G$. If $T$ is generic,
$\Phi(M^T)$ is a set of isolated points. In general, however, the
images $\Phi(M^{T})$ form {\em codimension-$k$ walls} of the polytope,
where $d-k$ is the dimension of the subtorus of $G/T$ acting
effectively. These walls may be internal to the polytope.  Let
$S^1\subset G$ be a 1-dimensional subtorus. If a component of
$M^{S^1}$ has an effective $G/S^1$ action, then its image will be a
codimension-1 wall. These include but are not restricted to the facets
of the polytope.

When $\mu$ is a regular value of $\Phi$, $\Phi^{-1}(\mu)$ is a
submanifold of $M$ and has a locally free $G$ action
by the invariance of $\Phi$. The quotient space $\Phi^{-1}(\mu)/G$ is called
the {\em symplectic reduction} and
is denoted
$M/\!/G (\mu)$, where the parameter $\mu$ is suppressed when
$\mu=0$.

This paper is concerned with the ordinary, rational cohomology of symplectic
reductions and its relationship to the $G$-equivariant cohomology of
$M$. By definition, the $G$-equivariant cohomology of $M$ is
$$
H_G^*(M) := H^*(M\times_G EG)
$$
where $EG$ is a contractible space with a free $G$ action and
$M\times_G EG$ indicates the product $M\times EG$ quotiented by the
diagonal $G$ action. For more details, see \cite{AB:momentmap}.

When $G$ acts locally freely on a manifold Z, the $G$-equivariant
cohomology of Z is the ordinary cohomology of the quotient $Z/G$. In
particular for $\mu$ a regular value of the moment map,
$$H_G^*(\Phi^{-1}(\mu)) = H^*(M/\!/G (\mu)).$$
\begin{theorem}[Kirwan]\label{th:surjective}
Let $M$ be a compact symplectic manifold with a Hamiltonian $G$
action, where $G$ is a compact torus. Let $\Phi$ be a moment map for
the $G$ action on $M$. For any regular value $\mu\in \g^*$, the
natural map
$$\kappa: H_G^*(M) \longrightarrow H^*(M/\!/G(\mu))$$
induced from the  inclusion $\Phi^{-1}(\mu)\subset M$ is a surjection.
\end{theorem}
We generalize this theorem to the following equivariant version.
\begin{theorem}\label{th:equvtsurjective}
Let $M$ be a compact symplectic manifold with a Hamiltonian $G$
action, where $G$ is a compact torus and $T$ a subtorus, not
necessarily generic. Let $\Phi_G$ and $\Phi_T$ be moment maps for the
respective actions.  For a regular value $\mu\in \mathfrak{t}^*$ of
$\Phi_T$, the inclusion of the submanifold $\Phi_T^{-1}(\mu)$ of $M$
induces a surjection in equivariant cohomology
$$
\kappa_T: H_G^*(M) \longrightarrow  H_{G/T}^*(M/\!/T
(\mu)).
$$
\end{theorem}
Note that we use $\kappa$ to indicate the map on cohomology induced by
the restriction from $M$ to $\Phi^{-1}_G(0)$, and $\kappa_T$ for the
map on cohomology induced by the restriction from $M$ to
$\Phi^{-1}_T(0)\supseteq \Phi^{-1}_G(0)$.

These two theoreoms allow one to compute the cohomology of the
symplectic reduction by a torus action, as long as one can compute the
equivariant cohomology and the kernel of this map. In certain cases,
the equivariant cohomology has been completely described
\cite{GKM:eqcohom}, \cite{GH:Gtocircle}.
The kernel of the Kirwan map in
Theorem~\ref{th:surjective} has also been described in very general terms, which we present below as Theorem~\ref{th:tw}. These results
hinge on the key fact that for compact Hamiltonian torus spaces, the
equivariant cohomology injects into the equivariant cohomology of the
fixed point set under the natural restriction map \cite{KiBook:quotients}.
\begin{theorem}[Kirwan] Let $M$ be a compact Hamiltonian $G$
space for $G$ a compact torus. Let $\mathcal{C}$ be the collection of
connected components of the fixed point set $M^G$. Then
\begin{equation}\label{eq:injection}
H_G^*(M)\hookrightarrow H_G^*(M^G)=\bigoplus_{C\in \mathcal{C}}
H_G^*(pt)\otimes H^*(C)
\end{equation}
is an injection. For a class $\alpha\in H_G^*(M)$ and $C\in
\mathcal{C}$, we write $\alpha|_C$ to indicate the restriction to the
fixed component $C$.
\end{theorem}

The relationship between surjectivity onto the symplectic reduction and 
injectivity of the equivariant cohomology into that of the fixed point
set is manifest in the following description of $\ker\kappa$ from
Theorem~\ref{th:surjective}.

\begin{theorem}[Tolman-Weitsman]\label{th:tw}
Let $M$ be a compact Hamiltonian $G$-space, with moment map
$\Phi$. Let $\mu\in\g^*$ be a regular value of $\Phi$.
The kernel of the Kirwan map
$$
\kappa:H_G^*(M) \longrightarrow H^*(M/\!/G(\mu))
$$
is the ideal $\langle K_G^\g(\mu)\rangle$ in $H_G^*(M)$ generated by $K_G^\g(\mu) = \bigcup_{\xi\in\g} K_G^\xi(\mu)$
where
\begin{align*}
K_{G}^\xi(\mu) := &\{ \alpha\in H_G^*(M)| \alpha |_C = 0 \mbox{ for
all connected components } C\mbox{ of } M^G \\
&\mbox{ with }\langle \Phi_G(C),\xi\rangle >\langle
\mu,\xi\rangle\}.
\end{align*}
\end{theorem}

We prove an equivariant analogue of this theorem:

\begin{theorem}\label{th:kernel} Let $M$ be a compact Hamiltonian $G$
space with moment map $\Phi_G$. Let $T$ be any subtorus of $G$
and $\Phi_T$ the corresponding moment map for the $T$ action. For
$\mu$ be a regular value of $\Phi_T$, the kernel of the equivariant
Kirwan map
$$
\kappa_T:H_G^*(M)\longrightarrow H_{G/T}^*(M/\!/T(\mu))
$$
is the ideal $\langle K_G^\mathfrak{t}(\mu)\rangle$
generated by $K_G^\mathfrak{t}(\mu) = \bigcup_{\xi\in\t} K_G^\xi(\mu)$ where $K_G^\xi(\mu)$
is defined as in
Theorem~\ref{th:tw}.

\end{theorem}
The subindex $G$ on $K_G^\g(\mu)$ and $K_G^\mathfrak{t}(\mu)$
indicates that these generate ideals in $H_G^*(M)$, while the
superindices indicate the relevant set of vectors in the theorems. We
will use $\langle K_G^\mathfrak{t}(\mu)\rangle$ to indicate
the ideal in $H_G^*(M)$ generated by $K_G^\mathfrak{t}(\mu)$ and we
will suppress the parameter $\mu$ when $\mu=0$. Notice that the
difference between the kernel of $\kappa$ and that of $\kappa_T$ is
that the union in
Theorem~\ref{th:kernel} is taken only over $\xi\in\t$.
The significance is that
Theorem~\ref{th:tw} can be recovered by the successive
application of Theorem~\ref{th:kernel} to one-dimensional subtori of
$G$. In the case that $\mu=0$, for each $S^1\subset G$, the kernel is generated by $K^\xi_{S^1}$
and $K^{-\xi}_{S^1}$ for
a choice of generator $\xi\in \mathfrak{s}^1$. It follows that the kernel of
$\kappa$ is generated by classes $\alpha\in H_G^*(M)$ satisfying one
of $d$ conditions, where $d=\dim G$.

The final contribution of this article is to find a small set $\Xi$ of
$\xi\in \g$ such that $K_G^\g(\mu)$ in Theorem~\ref{th:tw} can be
replaced by $K_G^{\Xi}(\mu):= \bigcup_{\xi\in\Xi}K_G^\xi(\mu)$.
Let $\Xi\subset \g$ be the (finite) set of unit
vectors perpendicular to codimension-1 walls of the moment polytope.
For any $\xi\in \Xi$, the annihilator $\xi^\perp$ is a hyperplane in
$\g^*$ through 0 and parallel to a codimension-1 wall of the moment
polytope. These hyperplanes, shifted to pass through $\mu$, are the
key to the kernel of $\kappa$.

\begin{theorem} The kernel of $\kappa:H_G^*(M)\rightarrow
H^*(M/\!/G(\mu))$ is generated by classes $\alpha\in H_G^*(M)$ with the following
property. There exists an oriented hyperplane $H_\mu^\alpha$ through $\mu$ and parallel to a
codimension-1 wall of the moment polytope such that $\alpha$ restricts to 0 on the
set of all fixed points whose image under $\Phi$ lie to the positive side of
$H_\mu^\alpha$.
\end{theorem}

We prove this theorem in Section~\ref{se:wallstheorem}.

\section{Equivariant Morse theory}\label{se:equvtMorse}

First we state several basic facts about equivariant Morse theory (as
developed in \cite{AB:YangMills}). We then refine these ideas to
gain equivariant homotopy information that is
standard in the case that the  Morse function has only isolated
critical points and there is no group action.

Let $f$ be a smooth function on a compact manifold $M$ and
$C$ a connected component of the critical set
of $f$ on $M$. Choose a Riemannian metric on $M$. We say $C$ is a {\em non-degenerate critical
manifold} for $f$ if
\begin{enumerate}
\item $C\subset M$ is a submanifold of $M$ such that $df=0$
along $C$, and
\item The Hessian $H_Cf$ (the matrix of second derivatives of $f$) is
non-degenerate on the normal bundle $\nu C$ of $C$ in $M$.
\end{enumerate}
If every connected component of the critical set
is non-degenerate, we say $f$ is {\em
Morse-Bott}.  At every non-degenerate critical submanifold, we use the Riemannian metric to 
identify a neighborhood of the zero-section in the
normal bundle $\nu C$ with a tubular neighborhood of $C$ in $M$.
The Hessian
defines a splitting $\nu C = \nu^+C\oplus \nu^-C$ into positive and
negative normal bundles. The dimension of the fibres of $\nu^-C$ is
called the {\em Morse index} of $C$ and is denoted $\lambda_C$.

Now assume that $M$ has a $G$ action, where $G$ is a compact torus,
and that the metric and the function $f$ are invariant with respect to this action. The
splitting of $\nu C$ into positive and negative bundles is
equivariant. This setting mimics that in which $f$ is a generic
component of the moment map, and $C$ is a connected component of $M^G$.

The proofs of the following two lemma are close to identical to those
presented by Milnor \cite{MiBook:Morsetheory} in the case that
critical manifolds are points and there is no group action. We note
only the minor differences in this more general setting. To make the
lemmas true for Morse-Bott functions, all local calculations must
include extra coordinates along the critical manifolds. To make the
lemmas equivariant, we choose an invariant Riemannian metric (so that the group acts by isometries)
and we equivariantly identify a neighborhood of the zero-section in the normal 
bundle of
a critical submanifold with a tubular neighborhood of that
submanifold (see \cite{AuBook}).
 This makes all relevant maps and homotopies
equivariant.
For Lemmas~\ref{le:Morselemma1} and \ref{le:Morselemma2},
assume $M$ is a compact
manifold with a $G$ action and a $G$-invariant Riemannian metric, where $G$ is a compact torus.
Let $f: M\rightarrow \R$ be an invariant Morse-Bott function.  Let $M^a = f^{-1}(-\infty,
a]$ for $a\in\R$.

\begin{lemma}\label{le:Morselemma1}
Suppose $a<b$ and $f^{-1}[a,b]$ is compact and contains no critical
points of $f$. Then
$$
M^a\simeq_G M^b
$$
is a $G$-equivariant homotopy equivalence.
\end{lemma}

\begin{lemma}\label{le:Morselemma2}
Suppose there are $m$ connected components $C_1,\dots C_m$ of the critical set
of $f$ with the same critical value, and such that
$f(C_1)=\cdots=f(C_m)\in(a,b)$. Then $M^b$ is equivariantly homotopic to $M^a$ with a
$\lambda_{i}$-cell bundle over $C_i$ attached for each $i=1,\dots m$, where $\lambda_{i}$ is the
Morse index of $C_i$.
\end{lemma}

These homotopy theorems lead the following results about
equivariant cohomology. For the sake of simplicity, we assume that
$f$ has distinct values for distinct connected components of the critical set. 
We note that the theorems 
proven below can be
generalized (using Lemma~\ref{le:Morselemma2}) to include the case that $f$ 
does not have distinct values for distinct critical sets, however the notation 
becomes cumbersome.
 
 Assume now that the $G$ action is Hamiltonian,
and that $f$ is a component of
the moment map for the $G$ action. As mentioned, we assume also that $f$
separates the critical sets, i.e.  one can order the critical sets $C_0,
C_1,\dots, C_k$ of $f$ so that $f(C_i)<f(C_j)$ if and only if $i<j$. If $f$ is generic, these
critical sets are the fixed points of the $G$ action.

The fundamental principle introduced by Atiyah and Bott
\cite{AB:YangMills} is that
\begin{theorem}
An equivariant cohomology class on $C_0$
extends to a class on $M^a$ for any $a\in \R$; in particular, it
extends to one on all of $M$ (although not uniquely).
\end{theorem}
We prove this theorem following \cite{AB:YangMills} and
\cite{TW:abelianquotient}.
\begin{proof}
Define the sets
\begin{align}\label{def:Mi+}
M_i^+ &:= f^{-1}(-\infty, f(C_i)+\epsilon_i)\mbox{ and}\\
\label{def:Mi-}
M_i^- &:= f^{-1}(-\infty, f(C_i)-\epsilon_i)
\end{align}
where $\epsilon_i>0$ is small enough that $C_i$ is the only critical
set in $f^{-1}(f(C_i)-\epsilon_i,f(C_i)+\epsilon_i)$.
Note that by Lemma~\ref{le:Morselemma1}
\begin{equation}\label{eq:homotopyequivalent}
M^+_i\simeq_G M^-_{i+1}.
\end{equation}
For each $i$, there is a long exact sequence in $G$-equivariant
cohomology
\begin{equation}\label{seq:relexact}
\cdots\rightarrow H_G^*(M^+_i,M^-_i)\rightarrow H_G^*(M^+_i)\rightarrow
H_G^*(M^-_i)\rightarrow H_G^{*+1}(M^+_i,M^-_i)\rightarrow\cdots.
\end{equation}
As before, we equivariantly identify a tubular neighborhood of $C_i$ with a
neighborhood of 0 in the normal bundle of $C_i$.
This bundle splits $\nu C_i = \nu^+C_i\oplus
\nu^-C_i$ into the positive and negative normal bundles of
$C_i$.  By excision and
homotopy equivalence,
$$
H_G^*(M^+_i,M^-_i) \iso H_G^*(D_i,S_i)
$$
where $D_i$ and $S_i$ are the unit disk and sphere bundles, respectively,
of $\nu^-C_i$. The equivariant Thom isomorphism states that
$$
H_G^*(D_i,S_i) \iso H_G^{*-\lambda_i}(C_i)
$$
where $\lambda_i = \dim(\nu^-C_i)$.
Thus the exact sequence (\ref{seq:relexact}) is equivalently
\begin{equation}\label{seq:Cilongexact}
\cdots\rightarrow H_G^{*-\lambda_i}(C_i)\rightarrow H_G^*(M^+_i)\rightarrow
H_G^*(M^-_i)\rightarrow H_G^{*+1-\lambda_i}(C_i)\rightarrow\cdots.
\end{equation}

\begin{lemma}\label{le:Cishortexact}
For $i=1,\dots, k$, the sequence (\ref{seq:Cilongexact}) splits into the
short exact sequence
\begin{equation}\label{seq:Cishortexact}
0\rightarrow H_G^{*-\lambda_i}(C_i)\rightarrow H_G^*(M^+_i)\rightarrow
H_G^*(M^-_i)\rightarrow 0.
\end{equation}
\end{lemma}
\begin{proof}[Proof of Lemma~\ref{le:Cishortexact}]
The composition $H_G^{*-\lambda_i}(C_i)\rightarrow
H_G^*(M^+_i)\rightarrow H_G^*(C_i)$, where the second map is induced by
inclusion, restricts to the composition
$H_G^{*-\lambda_i}(C_i)\iso H_G^*(D_i,S_i)\rightarrow
H_G^*(D_i)\iso H_G^*(C_i)$. This latter composition is muliplication
by the equivariant Euler class of the negative normal bundle of
$C_i$. Atiyah and Bott show in \cite{AB:YangMills} that, in the case
that there is an
$S^1\subset G$ which is acting on $\nu^-C_i$ and fixing $C_i$, this
class is not a zero-divisor. It follows that
$H_G^{*-\lambda_i}(C_i)\rightarrow H_G^*(M^+_i)$ must be an
injection. Furthermore, by the exactness of sequence
(\ref{seq:Cilongexact}), the image of $H_G^*(M^-_i)\rightarrow
H_G^{*+1-\lambda_i}(C_i)$ is the kernel of
$H_G^{*+1-\lambda_i}(C_i)\rightarrow H_G^{*+1}(M^+_i)$, which is 0 by
injectivity. Thus there is a surjection
$H_G^*(M^+_i)\rightarrow H_G^*(M^-_i),$ showing that the sequence splits.
\end{proof}
An equivariant cohomology class on $C_0$ extends to one on $M_0^+$. By
the homotopy equivalence (\ref{eq:homotopyequivalent}), a class on
$M_i^+$ extends to one on $M^-_{i+1}$. Surjectivity in
(\ref{seq:Cishortexact}) implies that a class on $M^-_{i+1}$ extends
to one on $M^+_{i+1}$. Thus by induction a class on $C_0$ extends to a
class in $H_G^*(M)$.
\end{proof}

One may ask the question of how unique these extensions are. By the injection
($\ref{eq:injection}$), a class is distinguished by its restriction to
the fixed point set. As these fixed point sets are critical sets for
Morse-Bott functions obtained from components of the moment map, we
exploit their relationship among each other.

Let $grad\ f$ be the gradient of $f$ with respect to a compatible Riemannian metric. For any critical subset
$C$, there is a cell-bundle of points $x\in M$
which converge to $C$ under the flow of $-grad\ f$ (or $grad\ f$), called the  
{\em stable manifold} (or {\em unstable manifold}) of $C$. Furthermore,
every point in $M$ converges to some $C$ under this flow. Thus, for
any $x\in M$, there is a (nonunique) sequence of critical sets $C_{i_1},
C_{i_2},\dots, C_{i_m}$ such that $x$ converges to $C_{i_1}$, and there
are points in the unstable manifold of $C_{i_j}$ which converge to $C_{i_{j+1}}$ for every
$j\geq 1$. Define the {\em extended stable manifold} of a critical set $C$ to be the set of
points $x\in M$ whose flows along $-grad\ f$  have an associated
sequence including $C$. In the case that $M$ is a coadjoint orbit of a semi-simple Lie group,
the extended stable manifold of a critical point $p$ is just the closure of the stable manifold
 out of $p$; they are called (permuted) Schubert varieties.

\begin{lemma}\label{le:restrictEuler}
Suppose $\alpha\in H_G^*(M)$ restricts to 0 on all $C_i$ such that
$i<j$. Then $\alpha|_{C_j}$ is some multiple of $e(\nu^-(C_j))$, the
equivariant Euler class of the negative normal bundle of $C_j$.
\end{lemma}

\begin{proof}
By (\ref{seq:Cishortexact}), a class $\alpha$ such that
$\alpha|_{M^-_j}=0$ is in the image of $H_G^{*-\lambda_j}(C_j)\rightarrow H_G^*(M_j^+)$.
As the map
$H_G^{*-\lambda_j}(C_j)\rightarrow H_G^*(M_j^+)\rightarrow H_G^*(C_j)$
(where the latter map is restriction) is multiplication by
$e(\nu^-C_j)$,  $\alpha|_{C_j}$ must be a multiple of this class.
\end{proof}

\begin{lemma}\label{le:flowupclass}
For every connected component $C$ of the critical set of $f$,
there is a class $\alpha$ with the following restriction properties:
\begin{enumerate}
\item\label{le:flowupclassrestriction1} $\alpha|_{C_i}=0$ if $C_i$ is not in the equivariant 
extended stable manifold of
$C$, and
\item\label{le:flowupclassrestriction2} $\alpha|_{C}=e(\nu^-C)$ where
$e(\nu^-C)$ is the equivariant Euler class of the
negative normal bundle (defined by $f$) of $C$.
\end{enumerate}
\end{lemma}

\begin{proof}
Let $j$ be the index such that $C_j=C$.
Note that $C_0,\dots, C_{j-1}$ are not in the extended stable manifold of $C_j$, as
$i<j$ implies $f(C_i)<f(C_j)$. Using the short exact sequence
(\ref{seq:Cishortexact}) we extend the class $0$ on $C_0$ to $\alpha$
on $M_j^-$ such that $\alpha$ restricts to 0 on $C_1,\dots,
C_{j-1}$. As $1\in H_G^{*-\lambda_j}(C_j)$ and the composition $H_G^{*-\lambda_j}(C_j)\rightarrow H_G^*(M_j^+)\rightarrow H_G^*(C_j)$
is multiplication by
$e(\nu^-C_j)$, we may extend $\alpha$ to $M_j^+$ such that its
restriction to $C_j$ is $e(\nu^-C_j)$.

Now suppose $i>j$ but $C_i$ is not in the extended stable manifold from
$C_j$.  Since points in the unstable manifold of $C_i$
flow into points in $M^-_j$, by Lemma~\ref{le:Morselemma2} there is a CW-complex $K$ which is
$G$-homotopic to $M_j^-$ with a $\lambda_{C_i}$-cell bundle over $C_i$
attached. $K$ has critical sets $C_1,\dots, C_{j-1}, C_i$. By
Lemma~\ref{le:restrictEuler} $\alpha|_K$ further restricts to
$m_ie(\nu^-C_i)$ on $C_i$ where $m_i\in H_G^*(C_i)$. Let $\beta_i\in
H_G^*(M)$ be such that $\beta_i|_{C_l}=0, l<i$ and
$\beta_i|_{C_i}=m_ie(\nu^-C_i)$. $\beta_i$ exists by the first part of
this proof, and $\alpha-\beta_i\in H_G^*(M)$ restricts to 0 on
$C_1,\dots, C_{j-1}, C_i$ and to $e(\nu^-C_j)$ on $C_j$. Let
$$
\gamma = \alpha-\sum\beta_i
$$
where the sum is over all $i>j$ such that $C_i$ is not in the extended stable manifold
of $C_j$. Then $\gamma$ has the desired restriction properties.
\end{proof}
Instead of dealing with CW-complexes, one might consider the closures
of the stable manifolds and make the same argument as that above with these
varieties. However, the singularities require resolving to be sure
there are well-defined classes, restricting one to the complex case only.

\section{The equivariant Kirwan map}
In this section we use equivariant Morse theory to prove
Theorem~\ref{th:equvtsurjective}, that the equivariant Kirwan map is
surjective. This proof follows very closely Tolman and Weitsman's
rendition of Kirwan's result (Theorem~\ref{th:surjective}) in the case
where $G=S^1$. The main elements of this proof are found in
\cite{TW:abelianquotient}, modified to allow for a torus action
commuting with the $S^1$ action and for the possibility of non-generic $T\subset G$.

\begin{proof}[Proof of Theorem~\ref{th:equvtsurjective}]
Choose $S^1\subset G$ and let
$$
\Phi_{S^1}:M\rightarrow (\mathfrak{s}^1)^*
$$
be a moment map for the $S^1$ action. We first show that the restriction
$H_G^*(M)\rightarrow H_G^*(\Phi_{S^1}^{-1}(0))$ is surjective.

Let $\xi\in\g$ generate the $S^1$ action. Consider the function
$(\Phi^\xi)^2$ where
$$\Phi^\xi:=\langle\Phi,\xi\rangle:M\rightarrow\R.$$
For an appropriate choice of norm on $(\mathfrak{s}^1)^*$, we have
$||\Phi_{S^1}||^2 = (\Phi^\xi)^2$. The critical set
of $(\Phi^\xi)^2$ consists of the minimum $\Phi_{S^1}^{-1}(0)$ and the
critical sets of $\Phi^\xi$. By the
moment map condition (\ref{eq:momentmapcondition}), the
critical points of  $\Phi^\xi$ are
the fixed points of the $S^1$ action generated by $\xi$. As Kirwan
notes in \cite{KiBook:quotients}, the function $(\Phi^\xi)^2$ may not be Morse-Bott
function; there may be degenerate critical
sets. However, this occurs only at the minimum $(\Phi^\xi)^{-1}(0)$
where the short exact sequence (\ref{seq:Cishortexact}) holds trivially.

As before, without loss of generality we may suppose that $f=(\Phi^\xi)^2$ separates the critical
 set. We order them by $f(C_i)<f(C_j)$ if and only if $i<j$ and $C_0 = \Phi^{-1}_{S^1}(0)$.  
Let $M^+_i$ and $M^-_i$ be as in
(\ref{def:Mi+}), (\ref{def:Mi-}).
We use Lemma~\ref{le:Cishortexact} to
show that $H_G^*(M^+_i)\rightarrow
H_G^*(\Phi_{S^1}^{-1}(0))$ is surjective for each $i$.

We noted that $H^*_G(M^+_0)\rightarrow H_G^*(C_0)$ is an
isomorphism. Now assume that we have a surjection
$H_G^*(M^+_i)\rightarrow H_G^*(C_0)$ for all $i\leq k-1$.  By the
short exact sequence (\ref{seq:Cishortexact}), there is a surjection
$H_G^*(M_k^+)\rightarrow H_G^*(M_k^-).$ But $H_G^*(M_k^-)\iso
H_G^*(M_{k-1}^+)$ by the homotopy
equivalence~(\ref{eq:homotopyequivalent}) and the latter ring surjects
onto $H_G^*(C_0)$ by
assumption. Thus $H_G^*(M_k^+)\rightarrow H_G^*(C_0)$ is a
surjection. By induction there is a surjection $H_G^*(M)\rightarrow
H_G^*(\Phi^{-1}_{S^1}(0)).$  As $0$ is a regular value of
$\Phi_{S^1}$, $H_G^*(\Phi^{-1}_{S^1}(0))\iso H^*_{G/S^1}(M/\!/S^1)$
and thus
\begin{equation}\label{eq:one-dimred}
H_G^*(M)\rightarrow H^*_{G/S^1}(M/\!/S^1)
\end{equation}
is a surjection.

There is a residual Hamiltonian $G/S^1$ action on $M/\!/S^1$, which
allows us to apply this technique inductively. By reduction in stages,
for any commuting subgroups $H_1$ and $H_2$ of $G$,
$(M/\!/H_1)/\!/H_2 = M/\!/(H_1\times H_2)$. For any $T\subset G$,
choose a splitting $T=S^1\times\cdots\times S^1$. Successively apply
the surjection (\ref{eq:one-dimred}) to obtain a sequence of
surjections
$$
H_G^*(M)\rightarrow H^*_{G/S^1}(M/\!/S^1)\rightarrow H^*_{G/(S^1\times
S^1)}(M/\!/(S^1\times S^1))\rightarrow\cdots\rightarrow H^*_{G/T}(M/\!/T).
$$
\end{proof}

\section{The kernel of the equivariant Kirwan map}

Here we prove Theorem~\ref{th:kernel}, that $\ker\kappa_T = K_G^\t(\mu)$, where
$\kappa_T:H_G^*(M)\rightarrow H_{G/T}^*(M/\!/T(\mu))$ is the equivariant
Kirwan map.
We first show that $K_G^\t(\mu)\subset \ker \kappa_T$ by
directly restricting classes in $K_G^\t(\mu)$ to $\Phi^{-1}_T(\mu)$. We then show
the ideals are equal by a dimension count.
\begin{proof}[Proof of Theorem~\ref{th:kernel}]
Let $i:\t\hookrightarrow\g$ be the inclusion map of Lie algebras and
$\pi:\g^*\rightarrow\t^*$ the induced projection. Denote by
$\langle,\rangle_G$ and $\langle,\rangle_T$ the natural pairings
between $\g^*$ and $\g$ and between $\t^*$ and $\t$,
respectively. Then if $\Phi_G:M\rightarrow \g^*$ is a moment map for
the $G$ action, $\Phi_T=\pi\circ\Phi_G$ is a moment map for the
restricted $T$ action. Choose $\mu\in\t^*$ a regular value of
$\Phi_T$.  Let $\alpha\in K_G^\xi$ for some $\xi\in\t$.
By definition $\alpha|_C = 0$ for every connected component $C$ of
$M^G$ such that $\Phi_G^{i(\xi)}(C)>\langle\mu,\xi\rangle_T$. Let
$$
M_\xi^+(\mu) = \{m\in M | \Phi_G^{i(\xi)}(m)>\langle\mu,\xi\rangle_T\}.
$$
$M_\xi^+(\mu)$ is a maximal dimension open $G$-invariant submanifold of
$M$. Thus the restriction of the injection $H_G^*(M)\rightarrow
H_T^*(M^G)$ to $M_\xi^+(\mu)$ is an injection into the cohomology
of those
components $C$ of $M^G$ which lie in $M_\xi^+(\mu)$. Thus $\alpha|_C=0$
for all $C\subset M_\xi^+(\mu)^G$ implies
$\alpha|_{M_\xi^+(\mu)}=0$. In particular
\begin{equation}\label{eq:res}
\alpha|_{\{m\in M|\ \langle\Phi_G(m),i(\xi)\rangle_G=\langle\mu,\xi\rangle_T\}}=0.
\end{equation}
But
\begin{align*}
\langle\Phi_G(m),i(\xi)\rangle_G &= \langle\pi\circ\Phi_G(m),\xi\rangle_T\\
& =\langle\Phi_T(m),\xi\rangle_T.
\end{align*}
Therefore,
\begin{align*}
\{m\in M|\
 &\langle\Phi_G(m),i(\xi)\rangle_G=\langle\mu,\xi\rangle_T\}\\
&= \{m\in M|\
 \langle\Phi_T(m),\xi\rangle_T=\langle\mu,\xi\rangle_T\}\\
&\supseteq
\Phi_T^{-1}(\mu).
\end{align*}
Thus $\alpha|_{\Phi_T^{-1}(\mu)} =0$, or equivalently
$\alpha\in\ker\kappa_T$. It follows that any class in the ideal generated by
$K_G^\mathfrak{t} = \cup_{\xi\in\mathfrak{t}} K_G^\xi$ lies in
$\ker\kappa_T$.

To show that the inclusion $\langle K_G^\mathfrak{t}\rangle
\subseteq \ker \kappa_T$ is an equality, we prove that
\begin{equation}\label{eq:dim}
\dim \langle
K_G^\mathfrak{t}\rangle
= \dim \ker \kappa_T = \dim(H_{G/T}^*(pt)\otimes \langle
K_T^\mathfrak{t}\rangle)
\end{equation}
 as graded ideals. We prove this in the case $\mu=0$, as the more
 general case is identical but notationally more cumbersome.

As $M$ is a Hamiltonian $G$-space and $M/\!/T$ is a Hamiltonian
$G/T$-space, they are both equivariantly formal with respect to
their group actions. This implies that as graded vector spaces
$H_G^*(M)=H_{G/T}^*(pt)\otimes H_T^*(M)$ and
$H_{G/T}^*(M/\!/T) = H_{G/T}^*(pt)\otimes H^*(M/\!/T).$
Recall that $\kappa_T:H_G^*(M)\rightarrow H_{G/T}^*(M/\!/T)$ and
$\kappa:H_T^*(M)\rightarrow H^*(M/\!/T)$ are
surjective (Theorems~\ref{th:equvtsurjective} and \ref{th:surjective}).
Thus there is a graded equality
\begin{align*}
\dim\ker \kappa_T &= \dim H_G^*(M)-\dim H_{G/T}^*(M/\!/T)\\
              &= \dim(H_{G/T}^*(pt)\otimes
              H_T^*(M))-\dim(H_{G/T}^*(pt)\otimes H^*(M/\!/T))\\
          &= \dim(H_{G/T}^*(pt)\otimes
              \ker(\kappa:H_T^*(M)\rightarrow H^*(M/\!/T)))\\
          &= \dim(H_{G/T}^*(pt)\otimes \langle
              K_T^\mathfrak{t}\rangle).
\end{align*}
where the last equality follows from Theorem~\ref{th:tw}.

We now show that $\langle K_G^\mathfrak{t}\rangle$ has the
same dimension. In degree $k$ for $\xi\in\t$,
$$
\dim \langle K_G^\xi\rangle^k = \dim\langle\{\alpha\in
H_G^k(M)|\ \alpha|_{C_i}=0\ \forall i<j\mbox{, for any $j$ such that }\Phi_G^\xi(C_j)>0\}\rangle
$$
where $i<j$ if and only if $\Phi_G^\xi(C_i)< \Phi_G^\xi(C_j)$. Let $F:H_G^*(M)\rightarrow H_T^*(M)$ be the
surjective map which forgets the $G/T$ action. Because $M^T=M^G$ for
generic $T$, if $\alpha\in H_G^k(M)$ has the property that
$\alpha|_{C_i}=0\ \forall i<j$, then
$F(\alpha)\in H_T^k(M)$ has the property that $F(\alpha)|_{C_i}=0\
\forall i<j$.
Furthermore, for $\xi\in\t$, $\Phi_T^\xi=\Phi_G^\xi$ so that such classes are precisely
 those in $\langle K_T^\xi\rangle$.
As $\ker F = H_{G/T}^*(pt)$ in all but degree 0, we conclude that
\begin{align*}
\dim \langle\{\alpha&\in
H_G^k(M)|\ \alpha|_{C_i}=0\ \forall i<j,\mbox{ for any $j$ such that }\Phi_G^\xi(C_j)>0\}\rangle \\
&= \dim\langle\{\beta\in\sum_{l+m=k}
H_{G/T}^l(pt)\otimes H_T^m(M)|\ \beta|_{C_i}=0\ \forall i<j,\\
&\hspace{1in} \mbox{ for any $j$ such that }\Phi_G^\xi(C_j)>0\}\rangle\\
&= \dim^k (H_{G/T}^*(pt)\otimes \langle K_T^\xi\rangle)
\end{align*}
\end{proof}

\section{The walls of the moment polytope}\label{se:wallstheorem}

In this section, we prove an important refinement of
Theorem~\ref{th:kernel}. It states that the collection of $\langle
K_G^\xi(\mu)\rangle$ for a small number of $\xi\in\g$ are sufficient to
generate the kernel of the map $H_G^*(M)\rightarrow
H^*(M/\!/G(\mu))$. In particular, one can consider only such $\xi$
which are perpendicular to codimension-1 walls of the moment polytope.

\begin{figure}[h]
\centerline{
\psfig{figure=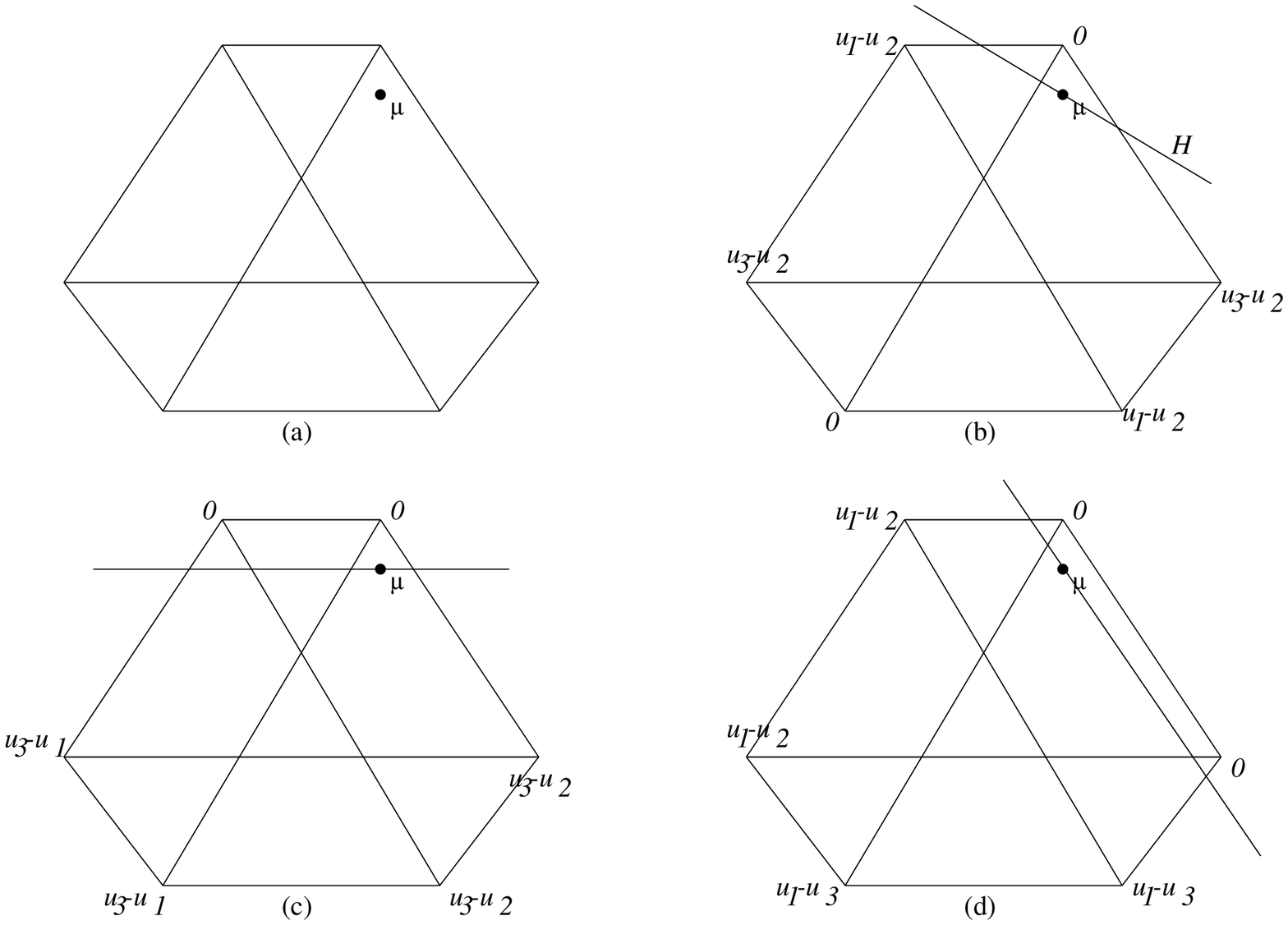,height=2.5in}
}
\centerline{
\parbox{4.5in}{\caption[]{\label{fig:SU3coadjointorbit}
{\small (a) The image of a generic $SU(3)$ coadjoint orbit under the
moment map; (b) A class $\alpha$ that is 0 to one side of a
codimension-1 hyperplane through $\mu$; (c)-(d) Two classes which
restrict to 0 on one side of a hyperplane parallel to a wall of the
moment polytope. Their sum is $\alpha$.}}}
}
\end{figure}

We illustrate the main theorem of this section with an example. Let
$M$ be a generic 6-dimensional coadjoint orbit of $SU(3)$. The maximal
torus $G=T^2$ of $SU(3)$ acts on this orbit in a Hamiltonian
fashion. The image of the moment map is a hexagon, and the
codimesion-1 walls of the moment polytope are shown in
Figure~\ref{fig:SU3coadjointorbit}(a).

Consider the reduction at the point $\mu$ indicated. According to
Theorem~\ref{th:tw}, a generating set for the kernel of the map
$H_G^*(M)\rightarrow H^*(M/\!/G(\mu))$ would include classes which
restrict to 0 on fixed points whose images under the moment map
(indicated by vertices) lie to one side of any hyperplane through
$\mu$. In particular, the class $\alpha$ whose restrictions to the fixed points
is indicated by Figure~\ref{fig:SU3coadjointorbit}(b) would be
a generator of the kernel because it is 0 to one side of the hyperplane
$H$. Theorem~\ref{th:wallssufficient} states that such a class is
redundant; it will in fact be generated by classes which are 0 to one
side of a hyperplane through $\mu$ {\em parallel to a wall} of the moment
polytope. In Figures~\ref{fig:SU3coadjointorbit}(c) and
\ref{fig:SU3coadjointorbit}(d) we see two classes which are in the
kernel according to Theorem~\ref{th:wallssufficient} and whose sum is
$\alpha$.

More generally, let $C$ be a connected fixed point component of the
$G$-action on $M$, and $\xi\in\g$ generic so that $C$ is a critical
manifold for $\Phi^\xi$. Let $X^\xi$ be the extended stable manifold of $C$
under $\Phi^\xi$.  The image $\Phi(X^\xi)$ of $X^\xi$ under the moment
map is not {\em a priori} convex. Suppose
$\dim\Phi(X^\xi)=\dim\g^*=d$. Consider the collection of codimension-1
walls of $\Phi(M)$ that lie in $\Phi(X^\xi)$. These are the images of
$S^1$-fixed point sets in $X^\xi$.  At least one of these walls,
translated to pass through $C$, has the property that $\Phi(X^\xi)$
lies entirely to one side. If $\Phi(X^\xi)$ is not maximal dimension,
then every codimension-1 wall containing $\Phi(X^\xi)$ has this
property.

We are now ready to state and prove in what sense the walls of the
moment polytope are sufficient information for calculating the kernel
of $\kappa$.

\begin{theorem}\label{th:wallssufficient}
Let $\Xi^\perp_\mu$ consist of hyperplanes through $\mu$ and parallel
to codimension-1 walls of the moment polytope. Let $K\subset H_G^*(M)$
be the ideal generated by classes $\alpha\in H_G^*(M)$ which restrict
to 0 on all connected components $C$ of $M^G$ whose images under
$\Phi$ lie to one side of some $H\in\Xi^\perp_\mu$. Then $K=\ker
\kappa$, where $\kappa:H_G^*(M)\rightarrow
H^*(M/\!/G(\mu))$ is the Kirwan map.
\end{theorem}

\begin{proof}

Choose any $\alpha\in H_G^*(M)$ where $\alpha|_{\Phi^{-1}(\mu)}=0$. By
Theorem~\ref{th:tw}, $\alpha\in\sum_{\xi\in\g}K^\xi_G(\mu)$. We want to show
that $\alpha$ can be written as a linear combination of elements in
$K^\eta_G(\mu)$ where $\eta\in\Xi$ are the annihilators of the hyperplanes
through $0$ and parallel to codimension-1 walls of the moment
polytope.

Without loss of generality, assume $\alpha\in K^\xi_G(\mu)$ for some $\xi$,
where $K^\xi_G(\mu)$ are classes restricting to 0 on fixed points whose image
under $\Phi$ lies to one side of $\xi^\perp+\mu$.  Order the connected
components of the critical sets $C_1,\dots, C_l$ so that $i<j$ if and
only if $\Phi^\xi(C_i)<\Phi^\xi(C_j)$.  We
prove that $\alpha$ can be expressed as a sum of elements in $K^\eta_G(\mu),
\eta\in \Xi$, by induction on the index of the critical sets.  Let
$C_{i_1}$ be the first critical set such that $\alpha|_{C_{i_1}}\neq
0$. Then $\alpha|_{C_{i_1}}$ is some multiple $m_{i_1}$ of
$e(\nu^-_{\Phi^\xi}C_{i_1})$.  Let $X_1^\xi$ be the extended stable manifold of
$C_{i_1}$ and let $\alpha_{i_1}$ be any class satisfying the
properties of Lemma~\ref{le:flowupclass}. In particular,
$\alpha_{i_1}\in K^{\eta_{i_1}}_G(\mu)\cup K^{-\eta_{i_1}}_G(\mu)$ where
$\eta_{i_1}$ is perpendicular to a codimension-1 wall $H$ of $\Phi(M)$
such that $\Phi(X_1^\xi)$
lies to one side of $H$ shifted to pass through $\mu$.  Then $\alpha-m_{i_1}\alpha_{i_1}$ is a class
which restricts to 0 on $C_1,\dots, C_{i_1}$. Now suppose that
$\alpha-\sum_{k=1}^{l}\alpha_{i_k}$ restricts to 0 on $C_1,\dots,
C_{i_l}$. Let $C_{i_{l+1}}$ be the first critical set on which
this class is non-zero. Use Lemma~\ref{le:flowupclass} to find a class
$\alpha_{i_{l+1}}$ supported on $X_{{l+1}}^\xi$. Then
$\alpha-\sum_{k=1}^{l+1}\alpha_{i_k}$ is 0 on $C_1,\dots,
C_{i_{l+1}}$. In this manner, we express $\alpha=\sum_i\alpha_i$,
where each $\alpha_i\in K^\eta_G(\mu)$ for some choice of $\eta$
perpendicular to a codimension-1 wall of the polytope.
\end{proof}

\section{Application to the product of symplectic manifolds}
Let $M=X_1\times\cdots\times X_k$ be the product of symplectic manifolds $X_i$,
each with a Hamiltonian $T$ action. Theorem~\ref{th:wallssufficient}
allows us to say a lot about the reduction of $M$ by the diagonal
torus action. We note that the diagonal torus $T_\Delta$ is a subtorus
of the product $G = T\times\cdots\times T$ acting on $M$. If the torus
is just one-dimensional, Theorem~\ref{th:wallssufficient} is not
more useful than the original formulation of Theorem~\ref{th:tw}. It
is when $\dim T \geq 2$ that one can significantly reduce the number
of vectors needed to generate the kernel of
$\kappa:H_{T_\Delta}^*(M)\rightarrow H^*(M/\!/T_{\Delta})$. We calculate the
cohomology of the reduced space for a product of two copies of
$\C P^2$s, where we have quotiented by the diagonal $T_\Delta^2$ action.

\begin{figure}[h]
\centerline{
\psfig{figure=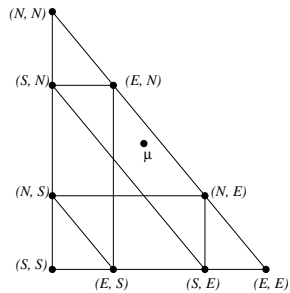,height=1.5in}
}
\centerline{
\parbox{4.5in}{\caption[]{\label{fig:cp2productmomentmap}
{\small The image of
the moment map for $T^2_\Delta$ acting on $\C P^2\times \C P^2$.}}}
}
\end{figure}

By the moment map condition (\ref{eq:momentmapcondition}),
if we choose the symplectic forms $\omega$ on $X_1$ and $k\omega$,
$k\in \R$ on $X_2$, the image for the moment map for the $T$ action on
$X_2$ is that for the action of $T$ on $X_1$ dilated by $k$. The image
of the moment map for the diagonal $T$ action on $X_1\times X_2$ is
the sum of the moment maps for each component. Let $X_1=X_2=\C P^2$
and $T$ act on $X_i$ by $(\theta_1,\theta_2)\cdot
[z_0:z_1:z_2]=[z_0: e^{i\theta_1}z_1: e^{i\theta_2}z_2]$.  Choose $k>2$ and let
$\Phi:X_1\times X_2\rightarrow \t^*_\Delta$ be a moment map for the diagonal
action. Then $\Phi(X_1\times X_2)$ with its walls is pictured in
Figure~\ref{fig:cp2productmomentmap}.

The image of fixed points are labeled and indicated by vertex dots, and the image
of fixed point sets of codimension-1 tori are indicated by line
segments, which in this case form codimension-1 walls of the moment polytope. We
find the cohomology of the symplectic reduction
$(X_1\times X_2)/\!/T_\Delta(\mu)$ where $\mu$ is the point indicated in the
figure by finding the kernel of $\kappa:H_{T_\Delta}^*(X_1\times
X_2)\rightarrow H^*((X_1\times X_2)/\!/T_\Delta(\mu))$ using
Theorem~\ref{th:wallssufficient}. We note that there are three distinct
hyperplanes through $\mu$ parallel to walls of the moment polytope:
horizontal, vertical, and diagonal with slope -1.

\begin{figure}[h]
\centerline{
\psfig{figure=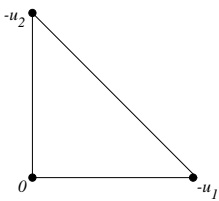,height=1in}
}
\centerline{
\parbox{4.5in}{\caption[]{\label{fig:symplres}
{\small The restriction of the class $x$ to the fixed point set of $T$
acting on $\C P^2$.}}}
}
\end{figure}

Consider first $H_{T_\Delta}^*(X_1\times X_2)$. By the equivariant K\"unneth
theorem,
\begin{equation}\label{eq:Kunneth}
H_{T_\Delta}^*(X_1\times X_2) = H_T^*(X_1)\otimes_{H_T^*(pt)}H_T^*(X_2),
\end{equation}
 where $T_\Delta\subset
G=T\times T$ is the diagonal torus. $H_T^*(\C P^2)$ is generated in
degree 2, by characteristic classes inherited from the module
structure $H_T^2(pt)\rightarrow H_T^2(\C P^2)$, and the equivaraint
symplectic form on $\C P^2$. We use the chosen basis for $\t^*\iso
H_T^2(pt)$ and denote the characteristic classes by $u_1$ and
$u_2$. These classes restrict to themselves on each fixed point of $T$
on $\C P^2$. Let $x$ be (the multiple of) the equivariant symplectic
class given by the restrictions indicated in
Figure~\ref{fig:symplres}.

Similarly, $u_1$ and $u_2$ are degree 2
equivariant classes on $\C P^2\times \C P^2$ which restrict to
themselves at every fixed point of the $T_\Delta$ action on the
product. By (\ref{eq:Kunneth}) the classes $u_1, u_2, x\otimes 1$ and
$1\otimes x$ generate $H_{T_\Delta}^*(\C P^2\times \C P^2)$. Here $u_i=u_i\otimes
1=1\otimes u_i$.

Note that a fixed point of $\C P^2\times \C P^2$ is a
pair $(p,q)\in ((\C P^2)^T, (\C P^2)^T)$. The restriction of a class
on the product space to a fixed point is
\begin{equation}\label{eq:productrestriction}
(a\otimes
b)|_{(p,q)}=a|_p\otimes b|_q
\end{equation} where $a,b\in H_T^*(\C P^2)$. The
algebraic structure of these classes is inherited  by multiplication
on each fixed point.

\begin{figure}[h]
\centerline{
\psfig{figure=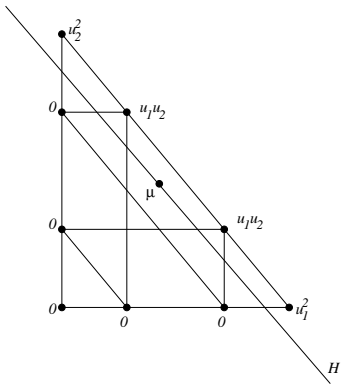,height=2in}
}
\centerline{
\parbox{4.5in}{\caption[]{\label{fig:x2restriction}
{\small The restriction of the class $x\otimes x$ to the fixed point
set of the diagonal $T$
acting on $\C P^2\times \C P^2$. Notice that this class restricts to 0
on fixed points whose image under $\Phi$ lie on one side of the
indicated hyperplane through $\mu$.}}}
}
\end{figure}

The Betti numbers for $H_{T_\Delta}^*(\C P^2\times \C P^2)$ are easy to
compute. As graded vector spaces, $H_{T_\Delta}^*(\C P^2\times \C P^2)=
H_T^*(pt)\otimes H^*(\C P^2\times \C P^2)$. It follows that the
equivariant Poincar\'e polynomial for $\C P^2\times \C P^2$ is
$(1+t^2+t^4+\dots)^2(1+t^2+t^4)^2 =
1+4t^2+10t^4+\dots.$ We noted above that $u_1, u_2, 1\otimes x,$ and
$x\otimes 1$ are four linearly independent degree 2 classes. A choice of 10
linearly independent degree 4 classes is $u_1^2, u_2^2,
u_1u_2, x\otimes u_1, x\otimes u_2, u_1\otimes x, u_2\otimes x,
x\otimes x, 1\otimes x^2,$ and $x^2\otimes 1$. The restrictions of
these classes to the fixed point set $(\C P^2\times \C P^2)^{T_\Delta}$ is
determined by the formula (\ref{eq:productrestriction}). As an
example, we compute the restriction of $x\otimes x$ to the fixed
point set. See Figure~\ref{fig:x2restriction}.

We have left to find the kernel of $\kappa:H_{T_\Delta}^*(\C P^2\times \C
P^2)\rightarrow H^*((\C P^2\times \C P^2)/\!/T_\Delta(\mu))$.
The reduced space $(\C P^2\times \C P^2)/\!/T_\Delta(\mu)$ has (real) dimension
$$\dim (\C P^2\times \C P^2) - 2\dim T_\Delta = 8-4=4.$$ Thus we expect all
classes of degree 6 or higher to be in the kernel of $\kappa: H_{T_\Delta}^*(\C
P^2\times \C P^2)\rightarrow H^*((\C P^2\times \C P^2)/\!/T_\Delta(\mu))$. In
degree 2, we can easily see that no linear combination of the four
classes above will restrict to 0 on one side of any of the three
hyperplanes parallel to walls of the moment polytope. By
Theorem~\ref{th:wallssufficient}, there are no degree 2 classes in
$\ker \kappa$. In degree 4, however, we expect nine (of ten) linearly
independent classes to be in the kernel, as the image of $\kappa$ is
(a multiple of) the volume form on the reduction. The reader can
verify using (\ref{eq:productrestriction}) that the following classes restrict to 0 on one side of one of
the three hyperplanes through $\mu$ parallel to a wall of the moment
polytope: $x\otimes x, u_1u_2+u_1\otimes x, u_1^2+u_1\otimes x,
u_1^2-1\otimes x^2, u_2^2-1\otimes x^2, u_2^2-u_2\otimes x, u_1\otimes
x-x\otimes u_1+u_2\otimes x-x\otimes u_2+x^2\otimes 1-1\otimes x^2,
1\otimes x^2+u_1\otimes x,$ and $x\otimes u_2+x\otimes u_1+x^2\otimes
1 + u_1u_2$.  They are all linearly independent, which can be verified
rather tediously. As an example, one can see in
Figure~\ref{fig:x2restriction} that $x\otimes x$ restricts
to 0 on one side of the diagonal hyperplane through $\mu$.


\end{document}